\documentclass{mymathart}
\usepackage[picturedir,theorems]{mymathmacros}
\usepackage{subfigure}

\numberwithin{equation}{section}

\usepackage{hyperref}

\renewcommand{\theta}{\vartheta}
\renewcommand{\phi}{\varphi}
\renewcommand{\rho}{\varrho}

\newcommand{\adds}{\underline{s}}

\newcommand{\s}{\adds}

\newcommand{\bdyit}[2]
             {{\rule{0pt}{0pt}_{\mbox{$\scriptstyle #2$}}^{\mbox{%
                   $\scriptstyle #1$}} }}

\renewcommand{\H}{\mathbb{H}}

\renewcommand{\P}{\mathcal{P}}

  \newcommand{\B}{\mathcal{B}}

\newcommand{\Blog}{\B_{\log}}

\newcommand{\V}{\mathcal{V}}

\title[On a question of Eremenko]{On a question of Eremenko \\
    concerning escaping components \\ of entire functions}

\author{Lasse Rempe}
\address{Dept. of Math. Sciences, University of Liverpool,Liverpool L69 7ZL, United Kingdom}
\email{lasse@maths.warwick.ac.uk}
\date{\today}
\thanks{This work was supported by the Research Development Fund at the
  University of Liverpool.}

\begin{document}

 \begin{abstract}
  Let $f$ be an entire function with a bounded set of singular values,
   and suppose furthermore that the postsingular set of $f$ is bounded. 
   We show that
   every component of the escaping set $I(f)$ is unbounded. This provides
   a partial answer to a question of Eremenko. 
 \end{abstract}

  \maketitle%

 \section{Introduction} \label{sec:introduction}
  
  Let $f:\C\to\C$ be a transcendental entire function, and consider the
   \emph{escaping set}
    \[ I(f) := \{z\in\C: f^n(z)\to\infty\} \]
   In 1989, Eremenko \cite{alexescaping} showed that $I(f)$ is nonempty,
    and that every component of $\cl{I(f)}$ is unbounded. He also
    states that the following seem plausible:
   \begin{enumerate} 
    \item Every component of $I(f)$ is unbounded. \label{item:eremenko}
    \item Every point of $I(f)$ can be connected to $\infty$ by
      a curve in $I(f)$. \label{item:strongeremenko}
   \end{enumerate}

  An interesting class of functions to consider for these questions
   is the \emph{Eremenko-Lyubich class}
    \[ \B := \{f:\C\to\C \text{ transcendental, entire}: 
                  \sing(f^{-1})\text{ is bounded}\}. \]
  (Recall that $\sing(f^{-1})$ consists of all critical and asymptotic
   values of $f$.)
   Recently \cite{fatoueremenko}, it was shown that 
   (\ref{item:strongeremenko}) fails, even in this class. In fact, 
   there exists $f\in\B$ for which $J(f)\supset I(f)$ contains no 
   nontrivial curve at all. On the other hand \cite{boettcher}, the
   dynamics ``near $\infty$'' remains the same throughout any given
   parameter space in $\B$ (in the sense of
   Eremenko and Lyubich; see \cite[Section 3]{alexmisha} or
   \cite[Section 2]{boettcher}). Taken together, these results show that
   there exist parameter spaces in class $\B$ in which
   (\ref{item:strongeremenko}) fails for every map. 

  In this note, we show that this is not the case
   for property (\ref{item:eremenko}). Recall that the
   postsingular set is defined as
    \[ \P(f) := \cl{\bigcup_{j\geq 0} f^j(\sing(f^{-1}))}. \]
   We shall prove the following result (which applies, in particular,
   to hyperbolic and postsingularly finite maps). 

 \begin{thm}[Escaping components of postsingularly bounded maps] 
    \label{thm:main}
  Suppose that $f\in\B$ and $\P(f)$ is bounded.  Then
   every component of $I(f)$ is unbounded.
 \end{thm}
 \begin{remark}
  Note that, if $g\in\B$ is arbitrary, then the function
   $f := g(z)/K$ will have a bounded postsingular set
   for sufficiently large $K$. Hence every parameter space in
   class $\B$ contains some functions to which our theorem
   applies. 
 \end{remark}
 
 There have previously been many cases in which property 
  (\ref{item:eremenko}) has been established
  by in fact proving the stronger property 
  (\ref{item:strongeremenko}), beginning with
  the case of a hyperbolic exponential map, which was completely described
  by Devaney in the early 1980s (see e.g.\ \cite{devaneykrych});
  the most recent result of this type 
  \cite[Theorem 1.2]{fatoueremenko} does so for any 
  (not necessarily hyperbolic) function $f\in\B$ of finite order. (Compare
  the latter reference for a more comprehensive discussion of previous
  contributions.) 
  It is interesting to note that each of these results
  also established the following property, which is
  stronger than (\ref{item:eremenko}):
   \begin{enumerate}  \setcounter{enumi}{2}
    \item For every $z\in I(f)$, there is an unbounded, connected set
      $A\ni z$ such that $f^n|_A\to \infty$ uniformly.
     \label{item:uniformeremenko}
   \end{enumerate}

  The methods of \cite{fatoueremenko} can be
   used to show the existence of hyperbolic entire functions $f\in\B$
   for which this property fails, so our proof will need to
   establish (\ref{item:eremenko}) without any uniform properties. 

  There is an interesting case of functions \emph{outside} class $\B$
   where this has previously been done. Rippon and Stallard
   \cite{ripponstallard} prove that, for \emph{any} transcendental
   entire function
   $f$, the set $A(f)\subset I(f)$ introduced by Bergweiler and Hinkkanen
   \cite{walteraimo} has only unbounded components (and this set does,
   in fact, satisfy the analog of (\ref{item:uniformeremenko})). In the
   case of an entire function which has a multiply-connected wandering domain,
   they also show that $A(f)$ is connected; since $A(f)$ is
   dense in $I(f)\cap J(f)$, this easily implies that $I(f)$ is also connected.
   This argument does not require (\ref{item:uniformeremenko});
   in fact, it seems plausible that this property will fail for
   functions 
   of this kind. 

  It is now well-known that Rippon and Stallard's construction of unbounded
   connected subsets of $I(f)$ can be adapted to give more precise
   information for functions of class $\B$. We will use a result of
   this type (Proposition \ref{prop:ripponstallard}) as a starting point
   of our proof;
   however, for the functions we are considering, $I(f)$ is usually
   disconnected.
   Thus, we need to find another way to connect an arbitrary
   escaping point to an unbounded escaping component. We shall achieve this
   by using three additional ingredients: the combinatorial
   coding of escaping points in class $\B$ by \emph{external addresses}, 
   hyperbolic expansion properties of our map $f$, and (perhaps
   most importantly)
   separation properties of the plane.

 \subsection*{Background and Notation}
  We refer the reader to \cite{jackdynamicsthird} for background on
   holomorphic dynamics and an introduction to plane hyperbolic geometry. 
   We denote the complex plane by 
   $\C$; all closures will be taken in $\C$ unless stated otherwise. 
   We denote the right 
    half-plane by
    $\H := \{\re z >0\}$;  
   more generally, we write $\H_{R} := \{\re z > R\}$. We also denote
   the unit disk by $\D$.
   If $U\subset\C$ is a domain and omits at least two points, we 
   will denote hyperbolic distance in $U$ by $\dist_U$, while Euclidean
   distance is denoted, as usual, by $\dist$.


 \section{Preliminaries}

  Let us define a \emph{tract} to be any Jordan domain in which
   $\re z$ is unbounded from above and which is disjoint from
   its $2\pi i \Z$-translates. 
   As in \cite{boettcher}, we denote by 
   $\Blog$ the set of all functions 
   $F:\V\to H$, 
   where 
   \begin{enumerate} 
    \item $H$ is a $2\pi i$-periodic Jordan domain which contains 
     a right half-plane.  
    \item $\V$ is $2\pi i$-periodic and $\re z$ is bounded from below,
     but not from above, in 
     $\V$. 
    \item Each component $T$ of $\V$ is a tract. (We call these components
      the \emph{tracts of $F$}).  
    \item
      For each such $T$,
      $F_T := F|_T$ is a conformal 
      isomorphism between $T$ and $\H$ with $F_T(\infty)=\infty$.
    \item The components of $\V$ accumulate only at $\infty$; i.e.,
      if $z_n\in\V$ is a sequence of points all belonging to different
      components of $\V$, then $z_n\to\infty$. 
   \end{enumerate}

  Recall that any entire function $f\in\B$ has a 
   \emph{logarithmic transform} $F\in\Blog$, with
   $\exp\circ F = f\circ\exp$. Thus the study of the dynamical 
   behavior of $f$ near $\infty$ reduces to that of $F$. 
   
  Note that any $F\in\Blog$ extends continuously to 
   $\cl{\V}$ by Carath\'eodory's theorem. We define 
   the escaping set of $F$ as
    \[ I(F) := \{z\in \cl{\V}: F^n(z)\in\cl{\V} \text{ for all
                    $n$ and } \re F^n(z)\to\infty\}. \]
   It is not difficult to see that $\cl{I(F)}$ 
    does not separate the plane and has
    no interior. Indeed, 
    this is a simple application of Koebe's
    distortion theorem, analogous to
    the proof of \cite[Theorem 1]{alexmisha}.
   
  A function $F\in\Blog$ is called \emph{normalized} if $H=\H$ and
   $|F'(z)|\geq 2$ for all $z\in \V$. (By 
   \cite[Lemma 1]{alexmisha}, any map
   $F\in\Blog$ can be normalized using a suitable restriction and
   conjugacy by a translation.)
   Let $z_0\in I(F)$, and let
   $T_j$ be the tract of $F$ with $F^j(z_0)\in \cl{T_j}$. Then the
   sequence
    $\s=T_0 T_1 T_2 \dots$ is called the \emph{external address} of
    $z_0$. The shift map on external addresses is denoted by
    $\sigma(\s) = T_1 T_2 \dots$.

 \section{Proof of the Theorem}

 Let us begin by noting a simple topological fact. 

\begin{figure}
 \begin{center}
 {\subfigure[Proof of Lemma \ref{lem:separation}\label{fig:windingtract}]{%
  \resizebox{.42\textwidth}{!}{%
  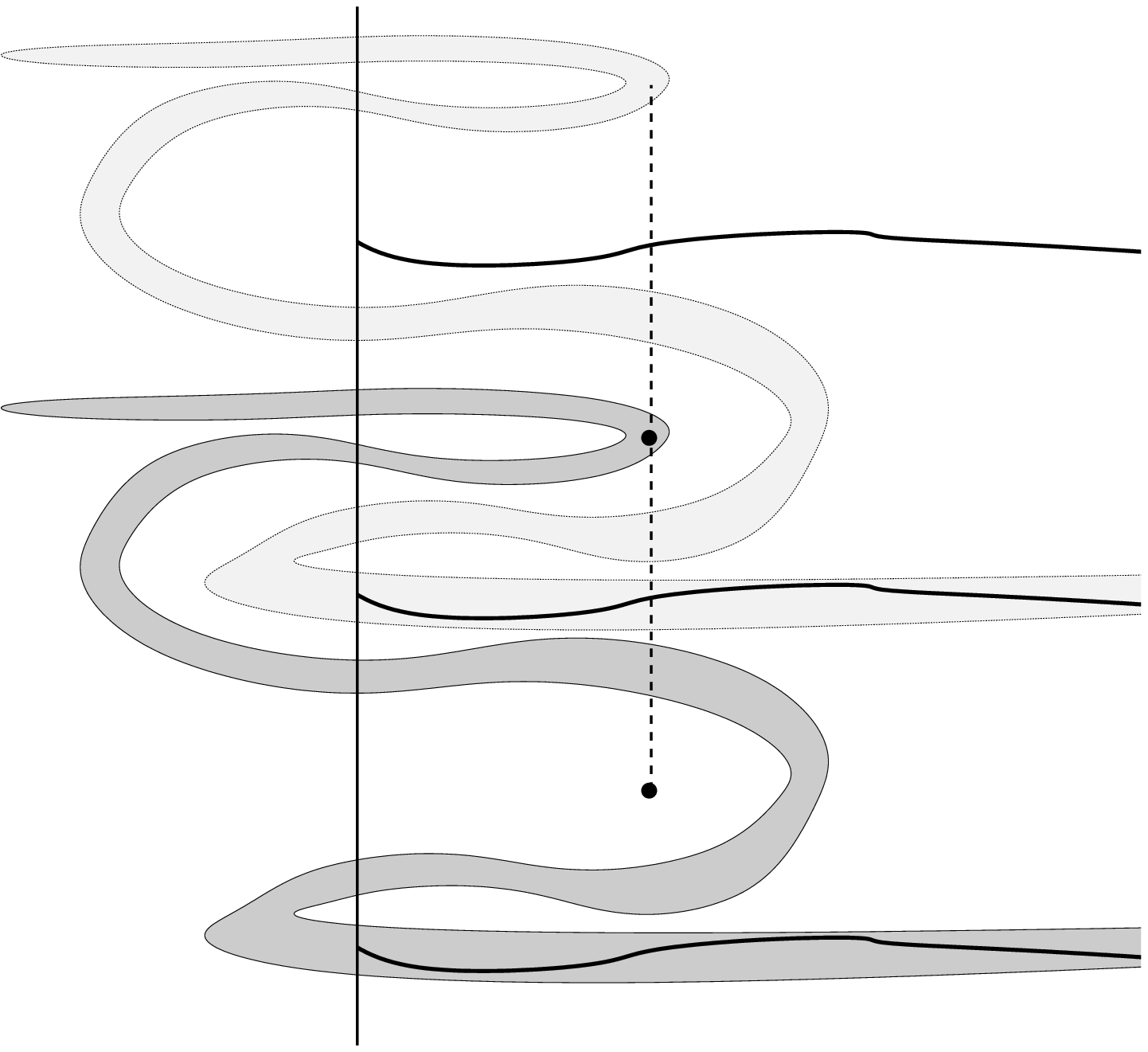}}}\hspace{1.5cm}
 {\subfigure[Construction of $A_j^k$ in 
           Proposition \ref{prop:ripponstallard}\label{fig:ripponstallard}]{%
  \resizebox{.42\textwidth}{!}{%
  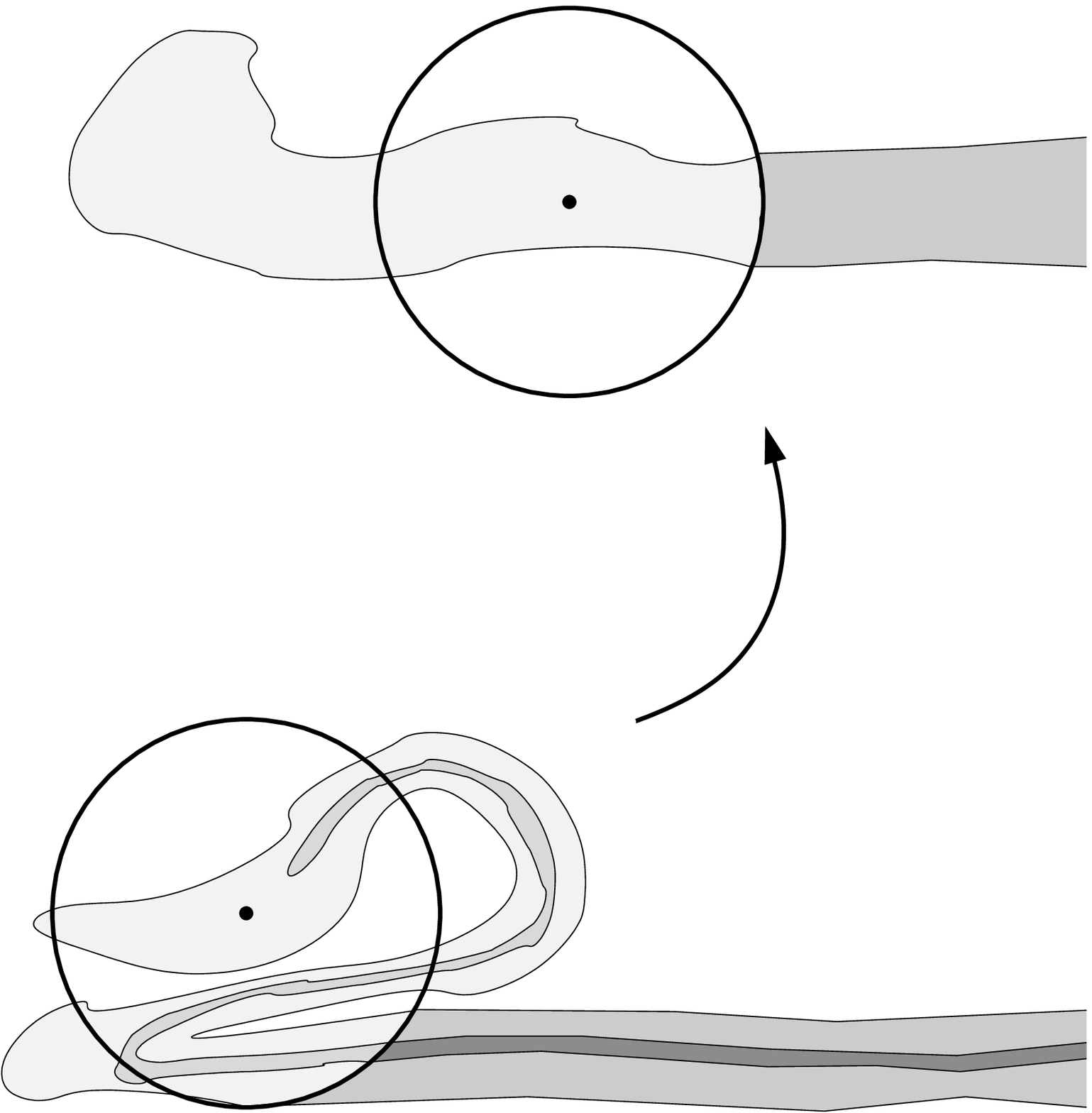}}}
  \end{center}
  \caption{Illustration of the proofs of Lemma \ref{lem:separation}
    and Proposition \ref{prop:ripponstallard}}
\end{figure}

 \begin{lem}[Separation of Tracts] \label{lem:separation}
  Let $T$ be a tract, let $R>0$ and suppose that $z$ belongs to
   the unbounded component $U$ of $\{z\in T: \re z > R\}$. Then
   the unbounded component
   $\wt{U}$ of 
   \[ T\setminus \{ z + it: t\in (-2\pi, 2\pi)\}. \] 
   is contained in $U$. 
 \end{lem}
 \begin{proof} We will prove the contrapositive. That is, suppose that
   $z\in T$ with
   $\re z> R$, and that there is 
   a curve $\gamma$ in
   $\wt{U}$ which connects $\infty$ to a point $w$ with $\re w = R$. 
   We may assume that all points of $\gamma$ except the endpoint $w$ 
   have real parts larger than $R$; then $\gamma$ belongs to $U$.
   We need to show that $z\notin U$.

  For some integer $j$,  the translates
   $\gamma_1 = \gamma + 2\pi i j$ and $\gamma_2 = \gamma + 2\pi i (j+1)$
   surround $z$ in $\H_{R}$
   (compare Figure \ref{fig:windingtract}). Clearly these translates must
   intersect the line segment $[z-2\pi i, z + 2\pi i]$. Since $\gamma\subset
   U$, it follows that $j\notin \{0,-1\}$. Hence $\gamma_1$ and
   $\gamma_2$ separate $z$ from $\gamma$, and hence from $U$, in
   $\H_R$. 
 \end{proof}

 Using this fact, we can provide an adaption of the
  aforementioned result of Rippon and Stallard \cite{ripponstallard}
  to class $\B$.
  (See also \cite[Theorem 3.3]{fatoueremenko}.)

\renewcommand{\labelenumi}{\theenumi}

 \begin{prop}[Unbounded sets of escaping points] \label{prop:ripponstallard}
  Let $F\in\Blog$ be normalized, and 
   let $z_0\in I(F)$ have external address
   $\s=T_0 T_1 T_2 \dots$. 

  Set 
   $D_j := \D_{2\pi}(F^j(z_0))$.
   Then
   there exists a sequence $(B_j)_{j\geq j_0}$ of
   unbounded subsets of $I(F)$ such that
  \begin{enumerate}
   \item[(i)] $B_j\subset T_j\setminus D_j$;
   \item[(ii)] $B_j\cup\{\infty\}$ is compact and connected;
   \item[(iii)] $B_j\cap \partial D_j\neq\emptyset$;
   \item[(iv)] $F(B_j)\subset B_{j+1}$. 
  \end{enumerate}
  \end{prop}
 \begin{proof}
   Let us set $z_j := F^j(z_0)$. Since $\re z_j\to\infty$, 
    we may assume without loss of generality that
    all $z_j$ belong to the unbounded connected component of
    $T_j\cap\H$.

   Let
   $A_j$ be the unbounded connected component of
   $\cl{T_j}\setminus D_k$. It follows from the previous lemma that
   $A_j\subset\H$, and furthermore 
   $\inf_{z\in A_j} \re z \to \infty$.

  We now define inductively $A_j^0 := A_j$ and let
   $A_j^{k+1}$ denote the unbounded component of
   $F_{T_j}^{-1}(A_{j+1}^k)\setminus D_j$. (In other words,
   $A_j^k$ is obtained by pulling back the set
   $A_{j+k}$ to $T_j$, cutting off at the disks $(D_i)$ in every step.)

  Then each $A_j^k$ is a closed, unbounded, connected
   set, and $A_j^{k+1}\subset A_j^k$. Furthermore, by the expanding
   property of $F$, we have $\dist(A_j^k,z_j)\leq 2\pi$. 

  We now set 
   \[ B_j := \bigcap_k A_j^k, \]
  and claim that these sets meet all our conditions. Indeed,
   $B_j \cup\{\infty\}$ is compact and connected as the
   nested intersection of compact and connected sets. Since each
   $A_j^k$ intersects $\partial D_j$, so does $B_j$. By construction,
   we have $B_j\subset \cl{T_j}\setminus D_j$, but in fact we also have
   $B_j\subset T_j$ since $F(B_j)\subset B_{j+1}\subset \H$. 
 \end{proof}


  The main problem in applying the previous proposition to prove our
   main result is that we do not know a priori that the 
   sets $B_j$ are connected. Hence it is conceivable that, for every
   component $C$ of $B_0$, there is some $j$ such that
   the component of $B_j$ containing $F^j(C)$ is far away from
   $z_j$. However, we will now show that the $B_j$ actually
   \emph{are} connected, so that this does not happen. The main tool
   required to do this is the following topological fact.

  \begin{lem}
   Let $T$ be a tract, and let $C_0,C_1\subset T$ be unbounded
    connected  sets. Then there exists $k\in\{0,1\}$ such that
    every point of $C_k$ has distance at most $2\pi$ from $C_{1-k}$.
  \end{lem}
  \begin{proof} Suppose that there is $z_0\in C_0$ such that
   $\dist(z_0,C_1)>2\pi$ (otherwise, there is nothing to prove). 
   Since $T$ is disjoint from its own $2\pi i\Z$-translates,
    we see that $C_1$ is disjoint from the set
    \[ A := [z_0 - 2\pi i , z_0 + 2\pi i] \cup \cl{C_0 + 2\pi i} \cup
         \cl{C_0 - 2\pi i}. \]

   Since $T$ is a Jordan domain which contains $C_0$ but is
    disjoint from $C_0+2\pi i$ and $C_0-2\pi i$, there is a component
    $U$ of $\C\setminus A$ which contains all sufficiently large points
    of $T$. In particular, we will have $C_1\subset U$. Furthermore,
    by construction of $U$, for every point $z\in U$ there is
    $t\in (-2\pi,2\pi)$ with $z+ti\in C_0$.  This completes
    the proof. 
  \end{proof}

   \begin{cor}[Uniqueness of unbounded components]
    Let $F\in\Blog$ be normalized, and suppose that $C_0,C_1\subset I(F)$
     are connected and unbounded. Suppose furthermore that
     all points in $C_0\cup C_1$ have the same external address. Then
     $C_0\cup C_1$ is connected.
   \end{cor}
   \begin{proof} Let $z_k\in C_k$ be two arbitrary points. By the previous
    lemma, for every $j\geq 0$ there is some $k\in\{0,1\}$ such that
    $\dist(F^j(z_k),F^j(C_{1-k}))\leq 2\pi$. In particular, there is
    some $k$ for which this happens for infinitely many $j$. Since
    $|F'|\geq 2$, it follows that
      \[ \dist(z_k,C_{1-k}) \leq 2^{1-j}\pi \]
    infinitely often; hence $\dist(z_k,C_{1-k})=0$. So $C_0\cup C_1$ is
    connected, as required.  
   \end{proof}

   In particular, this implies that the sets $B_j$ from Proposition
    \ref{prop:ripponstallard} are connected, as desired. Collecting
    these results together, we obtain:

   \begin{thm}[Unbounded escaping components] \label{thm:unboundedcomponents}
    Let $z_0\in I(F)$ have external address $\s=T_0 T_1 \dots$. 
    Then the set
      \[ I_{\s}(F) := \{z\in I(F): \text{$z$ has address $\s$} \} \]
    has a unique unbounded component $C_{\s}$. Furthermore,
      \[ \dist(F^j(z_0) , C_{\sigma^j(\s)}) \leq 2\pi \]
    for all sufficiently large $j$. (In fact, it suffices to take
    $j$ so large that $F^i(z_0)$ belongs to the unbounded component of
    $T_i\cap\H$ for all $i\geq j$.) 
   \end{thm}

 \begin{proof}[Proof of Theorem \ref{thm:main}]
   Let us assume without loss of generality that $0\in \P(f)$.
    By conjugating $f$ with $z\mapsto Kz$ for sufficiently large
    $K$, we may furthermore assume that
    $\P(f)\subset \D$.
    We set $W:= \{|z|>1\}$ and  $\wt{\V} := f^{-1}(W)$.
    Then every component of $\wt{\V}$ is an unbounded Jordan domain;
    these components are called the tracts of $f$. 
    We also set $U := \C\setminus \P(f)$; then $W\subset U$ by definition,
    and also 
    $f^{-1}(U)\subset U$. 

  Now let $z_0\in I(f)$. Note that finitely many iterates of $z_0$ might
   not belong to a tract of $f$; however, we will nonetheless be able
   to assign an external address to $z_0$. Let us begin by 
   associating to $z_0$ a sequence
   $\wt{T}_1, \wt{T}_2, \dots$ of tracts of $f$. 
   For sufficiently large $k$
   (say $k\geq k_0$), the point
   $f^k(z_0)$ will belong to the unbounded component of
   $\wt{T}\cap W$ for some tract $\wt{T}$.
   (Note that only finitely many tracts can intersect $\partial W$, and
    recall that $|f^n(z_0)|\to \infty$.) We will denote this tract by
    $\wt{T}_k$. To define $\wt{T}_j$ for $j<k_0$ is a little bit
    trickier. To do so, take $k\geq k_0$ and let $\wt{V}_k$  be the 
    component of $f^{-k}(\wt{T}_k)$ containing $z_0$. Then
    $\wt{V}_k$ is a Jordan domain, and 
    (since $\wt{T}_k$ is disjoint from the postsingular set)
    $f^k:\wt{V}_k\to \wt{T}_k$ is
    a conformal isomorphism. Furthermore, for $j< k$, all sufficiently
    large points of $f^j(\wt{V}_k)$ will be contained in a
    unique tract
    $\wt{T}_j$. This completes our definition of the sequence
    $(\wt{T}_j)$. (Note that the definition is independent of $k$, and that
    it agrees with our original definition for $j\geq k_0$.)

   However, a sequence of tracts of $f$ is not quite sufficient to 
    identify a suitable escaping component. 
    The easiest way to adapt this is to
    instead consider a logarithmic transform $F\in\Blog$. More
    precisely, let 
      $\V := \exp^{-1}(\wt{\V})$. Then there is a function
    $F:\V\to \H$ with 
    $\exp\circ F = f\circ \exp$. If $K$ was chosen sufficiently large
    in the beginning of the proof,
    then $F$ will be normalized; i.e.\
    $|F'| \geq 2$. (See \cite[Section 2]{alexmisha} or
    \cite[Section 2]{boettcher} for more details on this logarithmic
    change of variable.)

   We now define an external address $\s=T_1 T_2\dots$ of $F$
    with $\exp(T_j)=\wt{T}_j$ as follows. Let $T_1$ be an arbitrary
    component of $\exp^{-1}(\wt{T}_1)$, and let $V_k$ be the
    unique component of $\exp^{-1}(\wt{V}_k)$ which has an unbounded
    intersection with $T_1$. Then $T_j$ is defined to be the
    unique tract of $F$ with $F^j(\zeta)\in T_j$ for all sufficiently
    large $\zeta\in V_k$ ($k\geq j$). 

  To complete the proof, let 
   $C_{\s}$ be the unbounded component of $I(F)$ having address
   $\s$, and let $C$ be the component of $I(f)$
   containing $\exp(C_{\s})$. For $k\geq k_0$, let $\zeta_k$ be
   the unique point in $T_k$ with $\exp(\zeta_k)= f^k(z_0)$. 
   We have $\dist(\zeta_{k+1},C_{\sigma^{k+1}(\s)})\leq 2\pi$ by 
   Theorem \ref{thm:unboundedcomponents}. 
   Since $\re z_{k+1}\to+\infty$, this means that
   the hyperbolic distance $\dist_{\H}(\zeta_{k+1},C_{\sigma^{k+1}(\s)})$
   tends to zero. The map 
   \[ \exp\circ F_{T_k}^{-1} : \H\to \wt{T_k} \]
   is a conformal isomorphism, and it follows that
   \[\dist_{\wt{T_k}}(f^k(z_0),f^k(C))\to 0. \]

  Finally, because $f^k:\wt{V}_k\to \wt{T_k}$ is a conformal isomorphism,
   we have
   \[ \dist_U(z_0,C) \leq \dist_{T_k}(z_0,C) = \dist_{\wt{T_k}}(f^k(z_0),f^k(C))\to 0. \]
   Thus $z_0\in C$, as required. 
  \end{proof}

 Finally, let us note the following interesting corollary of Theorem
  \ref{thm:unboundedcomponents} for \emph{disjoint-type} maps.

 \begin{cor}[Escaping components of disjoint-type maps]
  Suppose that $F\in\Blog$ is  of
   \emph{disjoint type}
   (that is, $\cl{\V}\subset \H$). Then the set 
   $I_{\s}(F)$ is connected and unbounded for every external
   address $\s$.
 \end{cor}
 \begin{proof}
  This follows by a similar (but simpler) hyperbolic contraction argument
  to that in the previous proof. Indeed, in this setting we have
  $F(I_{\s})=I_{\sigma(\s)}$ (rather than the inclusion which usually holds).
  So if $z_0\in I_{\s}$, then 
  $\dist(F^n(z_0),F^n(I_{\s}))\leq 2\pi$, and $F$ is 
  again a (even strict)  hyperbolic expansion. The claim follows. 
 \end{proof}

 \begin{remark}
  \cite[Theorem 1.3]{boettcher} states that every \emph{hyperbolic} map
   $f\in\B$ is actually conjugate, on its escaping set, to a
   map of disjoint type. Hence the previous result
   provides a slightly 
   alternative proof of Theorem \ref{thm:main} in the hyperbolic case.
 \end{remark}

\section*{Acknowledgments}
  I would like to thank Walter Berg\-weiler, Helena Mi\-hal\-je\-vic-Brandt,
   J\"orn Peter, Phil Rippon and  
   Hendrik Schubert 
   for interesting discussions about this work.

\bibliographystyle{hamsplain}
\bibliography{../../Biblio/biblio}%

\end{document}